\documentclass[12pt]{amsart}

\def \fg{{\mathfrak g}}
\def \fh{{\mathfrak h}}
\def \fk{{\mathfrak k}}
\def \fp{{\mathfrak p}}

\def \fa{{\mathfrak a}}

\def \fn{{\mathfrak n}}
\def \fq{{\mathfrak q}}
\def \i{{\bf i}}

\newtheorem{theorem}{Theorem}[section]
\newtheorem{lemma}[theorem]{Lemma}
\newtheorem{proposition}[theorem]{Proposition}
\newtheorem{corollary}[theorem]{Corollary}
\newtheorem{definition}[theorem]{Definition}

\numberwithin{equation}{section}

\begin{document}

\baselineskip=16pt

\title[Lagrangian loci of quotients]{A note on Lagrangian loci of quotients}

\author[P. Foth]{Philip Foth}

\address{Department of Mathematics, University of Arizona, Tucson, AZ 85721-0089}

\email{foth@math.arizona.edu}



\date{March 21, 2003}

\begin{abstract}
We study hamiltonian actions of compact groups in the presence of 
compatible involutions. We show that the lagrangian fixed point set
on the symplectically reduced space is isomorphic to the disjoint
union of the involutively reduced spaces corresponding to 
involutions on the group strongly inner to the given one.   
Our techniques imply that the solution to the eigenvalues of a sum problem
for a given real form can be relayed to the quasi-split real form in the
same inner class. We also consider invariant quotients with respect to 
the corresponding real form of the complexified group. 
\end{abstract}

\maketitle

\section{Introduction}

Let a compact group $K$ act in a hamiltonian way on a manifold $X$. 
Let $\sigma$ be an anti-symplectic involution on $X$ and let $\tau$ 
be an involution on $K$ such that for any $x\in X$ and $k\in K$ we
have $\sigma(k.x)=\tau(k).\sigma(x)$. This setup was considered 
by O'Shea and Sjamaar in \cite{OSS} in order to establish a real 
analogue of Kirwan's convexity theorem which lead them to 
finding inequalities on possible spectrum of sums of two matrices. 

The main goal of the present paper is to establish a relationship
of the fixed point set on the reduced space $(X//K)^{\sigma}$ 
with the involutively reduced space $X^{\sigma}//K^{\tau}$. The 
latter is defined as the quotient by $K^{\tau}$ of the zero level 
set in $X^{\sigma}$. Our main observation is that one needs to 
consider simultaneously all conjugacy classes of involutions on $K$
which are {\it strongly inner} to $\tau$. We say that $\tau_s$ 
is strongly inner to $\tau$, if in addition to being inner to $\tau$, i.e. 
$\tau_s={\rm Ad}_s\circ\tau$, we require $s\tau(s)=1$. Then we show that
when the action of $K$ on the zero level set is free, the space
$(X//K)^{\sigma}$ is the disjoint union of such 
$X^{s\sigma}//K^{\tau_s}$, where $s$ runs through the connected components 
of the subset $Q$ of $K$ consisting of elements satisfying $\tau(k)=k^{-1}$. 
The elements of $Q$ correspond to involutions on the group $K$ which are 
strongly inner to $\tau$. We also discuss the singular case. 

We generate many examples by taking a complex semisimple  Lie group $G$ 
and its flag manifold $G/P$ with the property that the real 
dimension of the closed orbit of a real form $G^{\tau}$ on $G/P$ 
equals the complex dimension of $G/P$. Then we show that there exists
a symplectic structure on $G/P$ and an anti-symplectic involution
$\sigma$ on $G/P$ which is compatible with $\tau$. To get interesting 
examples of reduced spaces, one then might take a product of several 
complex flag manifolds of the above type with the diagonal $G$-action. 
A particular case when $G^{\tau}={\rm SL}(2,{\mathbb H})$ was considered
in \cite{FL}, where the involutively reduced space was identified 
with the moduli space of polygons in ${\mathbb R}^5$. 

As an application of our techniques we can show that the problem of 
finding a solution to $A_1+\cdots + A_k=0$, where $A_i$ is fixed by $\tau$ 
and has a prescribed spectrum, can be relayed to finding a solution
for a quasi-split involution in the same inner class. 

We can also consider a similar setup, when a linearly reductive 
complex Lie group $G$ equipped with an anti-holomorphic involution
acts on a projective variety $X$, which has a compatible anti-holomorphic
involution. We employ results of Richardson and Slodowy \cite{RS} on 
minimum vectors to establish an analogue of the aforemention result in
this context. 

In a subsequent work, we wish to apply our results in order to 
compute ${\mathbb Z}/2$-cohomology of the real loci of quotients, and  
extend the Kirwan surjectivity theorem for equivariant cohomology
of the real loci in the non-abelian case. For the case of a torus action, 
such a generalization was studied by Goldin and Holm \cite{GH}. 

\section{Involutions and reduction.}

Let $K$ be a compact connected Lie group and let $\tau: K\to K$ be a group 
homomorphism satisfying $\tau^2={\rm Id}$. We will refer to such a $\tau$ 
simply as an {\it involution} on $K$. Let $(X, \omega)$ be a hamiltonian 
$K$-manifold with an equivariant momentum map $\mu: X\to \fk^*$. Besides,
let $\sigma: X\to X$ be an involution on $X$, i.e. a diffeomorphism 
satisfying $\sigma^2={\rm Id}$. We say that $\sigma$ is an {\it anti-symplectic} 
involution if $\sigma^*\omega=-\omega$. We also shall say that 
$\tau$ and $\sigma$ are {\it compatible} if 
$$\sigma(k.x)=\tau(k).\sigma(x)\ \ .$$ 
We recall that \cite[Lemma 2.2]{OSS} asserts that $\mu$ can be shifted by 
an ${\rm Ad}^*$-invariant element of $\fk^*$ to ensure that 
$\mu(\sigma(x))=-\tau^*(\mu(x))$. We therefore shall assume that our $\mu$ 
satisfies this property, which in particular implies that $\mu^{-1}(0)$ 
is $\sigma$-stable.  

We denote by $K^{\tau}$ the set of points in $K$ fixed by $\tau$ and we will 
denote by $X^{\sigma}$ the set of points in $X$ fixed by $\sigma$. We notice that
$X^{\sigma}$ is a smooth submanifold of $X$ and that $K^{\tau}$ has a natural 
structure of a Lie subgroup of $K$. Recall that the Marsden-Weinstein
symplectic quotient $X//K$ of $X$ with respect to the $K$-action
is defined as $$X//K=\mu^{-1}(0)/K\ \ .$$ In the situation when $K$ acts freely 
on $\mu^{-1}(0)$, the value $0\in\fk^*$ is regular and the quotient $X//K$ 
has a natural structure of a smooth manifold carrying the reduced symplectic
form $\omega_{red}$. When $0$ is not a regular value of $\mu$, the symplectic 
quotient $X//K$ is a {\it stratified symplectic space} in the terminology of 
Sjamaar and Lerman \cite{SL}. It is easy to check that the involution $\sigma$ 
descends to an involution, also called $\sigma$, on $X//K$, where it is 
anti-symplectic with respect to $\omega_{red}$. Let us denote the fixed point 
set of the latter involution by $(X//K)^\sigma$. It is easy to see \cite{OSS}
that the fixed point set of an anti-symplectic involution is either empty or 
a lagrangian submanifold. In addition, $X^{\sigma}$ is $K^{\tau}$-stable and
its image in $\fk^*$ under the momentum map $\mu$ is contained in 
$(\fk^{\tau})^\perp$, which can be identified with $\fq^*$, where 
$\fq\subset \fk$ is the $(-1)$-eigenspace of $\tau$. 

\medskip\noindent{\bf Remark.} Our use of $\omega$ and $\mu$ is rather 
dictated by convenience than necessity. One can use a more general set-up,
where $\sigma$ is not required to be anti-symplectic and get same results. 
\medskip

Let us define the {\it involutively reduced space}:
$$X^{\sigma}//K^{\tau} = \left( \mu^{-1}(0)\cap X^{\sigma}\right) /K^{\tau}\ \  .$$
This definition makes sense due to the fact that 
$\mu^{-1}(0)\cap X^{\sigma}=\left( \mu^{-1}(0)\right)^\sigma$
is $K^{\tau}$-stable. 

Given any point in $X^{\sigma}//K^{\tau}$ 
represented by a $K^{\tau}$-orbit $K^{\tau}.x$ in $X$, we can consider 
the $K$-orbit $K.x$, which belongs to $\mu^{-1}(0)$ and thus defining a 
point in $X//K$. Moreover, since $\sigma(K.x)=K.x$,  
the corresponding point in $X//K$ is fixed by $\sigma$. The map thus obtained 
is denoted by: 
\begin{equation}
\psi: X^\sigma//K^{\tau}\to (X//K)^{\sigma}\ \ .
\label{Eqpsi}
\end{equation} 

\medskip
\noindent{\bf Remark.} Our involutively reduced space is related to the 
{\it Lagrangian quotients} defined in \cite{OSS}. Under the assumption that 
$\psi$ is injective, they are actually the same. 
\medskip

We wish to study properties of the map $\psi$. 

\begin{lemma} If $K$ acts freely on $\mu^{-1}(0)$, then $\psi$ is injective.  
\label{leminj} \end{lemma} 
\noindent{\bf Proof.} Assume that two $K^{\tau}$-orbits $K^{\tau}.x_1$ and
$K^{\tau}.x_2$ from $X^{\sigma}\cap\mu^{-1}(0)$ are mapped to the same
point by $\psi$. This would imply that there exists $k\in K$ such that 
$x_1=k.x_2$. Applying $\sigma$ to both sides: 
$$x_1=\sigma(x_1)=\tau(k).\sigma(x_2)=\tau(k).x_2 .$$ The fact that $K$ acts
freely on $\mu^{-1}(0)$ amounts to $k=\tau(k)$, meaning that
$x_1$ and $x_2$ are actually in the same $K^{\tau}$-orbit. 
\hfill {\bf Q.E.D.} \medskip

Denote by $Q$ the subset of elements in $K$ satisfying 
$\tau(k)=k^{-1}$. Let us denote by $Q_0$ the connected component of $Q$ 
containing the identity element. The group $K$ acts on $Q$ via 
$k.q=kq\tau(k)^{-1}$ and the connected components of $Q$ are precisely the 
orbits of this $K$-action. In particular, $Q_0$ is the collection of elements 
in $K$ representable as $k\tau(k)^{-1}$. It is easy to see that $Q_0$ 
is diffeomorphic to the symmetric space $K/K^{\tau}$. Actually, it is true 
that each connected component of $Q$ is also a symmetric space, for  
each connected component $Q_s$ of $Q$ and an
element $s\in Q_s$ we have:

\begin{lemma} For $s\in Q_s$, the map $\tau_s={\rm Ad}_s\circ\tau$
is an involution on $K$. 
\end{lemma}
\noindent{\bf Proof.} Due to $\tau(s)=s^{-1}$, we have
$$\tau_s\circ\tau_s(k)={\rm Ad}_s\circ\tau\circ{\rm Ad}_s\circ\tau(k)= 
{\rm Ad}_s(\tau(s)k\tau(s)^{-1})=k\ \ .$$
\hfill {\bf Q.E.D.} \medskip

Since any other element $s'\in Q_s$ can be presented as $ks\tau(k)^{-1}$, 
we see that the actions of $K$ on $Q_s$ and on $K/K^{\tau_s}$ are compatible.
Moreover, the connected component $Q_s$ is isomorphic to the symmetric space 
$K/K^{\tau_s}$. See \cite[Proposition 5.8]{Xu} for detailed discussion. 

\begin{definition} We say that an involution $\tau'$ on $K$ is {\em inner} 
to $\tau$ if there exists an element $s\in K$ such that 
$\tau'={\rm Ad}_s\circ\tau$. We also say that an involution $\tau'$ on $K$ is 
{\em conjugate} to $\tau$, if there exists an element $k\in K$ such that 
$\tau'={\rm Ad}_k\circ\tau\circ {\rm Ad}_{k^{-1}}$. 
\end{definition}

\begin{lemma} If $\tau'$ is conjugate to $\tau$, then it is also inner 
to it. If $\tau'={\rm Ad}_s\circ\tau$ is inner to $\tau$, then
$s\tau(s)$ belongs to the center of $K$. The properties of being 
conjugate or inner are equivalence relations.
\end{lemma}
\noindent{\bf Proof.} Straightforward.  \hfill {\bf Q.E.D.} \medskip

\begin{definition} We say that an involution $\tau'$ on $K$ is 
{\em strongly inner} to $\tau$ if the element $s$ in the definition can
be chosen from $Q$, i.e. $s\tau(s)=1$. 
\end{definition}

For example, if $K$ is of adjoint type, then being 
inner is the same as being strongly inner. 

\begin{lemma} Being strongly inner is an equivalence relation. 
\end{lemma} 
\noindent{\bf Proof.} Straightforward.  \hfill {\bf Q.E.D.} \medskip

Given $s\in Q$ and an involution $\tau_s={\rm Ad_s}\circ\tau$ strongly 
inner to $\tau$, we can introduce a diffeomorphism $\sigma_s$ on $X$:
$$\sigma_s(x)=s.\sigma(x)\ \ .$$

\begin{lemma} If $\tau_s$ is strongly inner to $\tau$, then
the map $\sigma_s$ is involutive, anti-symplectic, 
and compatible with $\tau_s$.  
\end{lemma}
\noindent{\bf Proof.} Involutivity:
$$ \sigma_s\circ\sigma_s(x)=s.\sigma(s.\sigma(x))=
s\tau(s).\sigma\circ\sigma(x) = x\ \ .$$
Anti-symplecticity readily follows from the facts that $\sigma$ is
anti-symplectic and $K$ acts in a Hamiltonian way. 
Compatibility: 
$$ \sigma_s(k.x)=s.\sigma(k.x)=s\tau(k).\sigma(x)=
{\rm Ad}_s\circ\tau(k)s.\sigma(x)=\tau_s(k).\sigma_s(x) .$$ 
\hfill {\bf Q.E.D.} \medskip

We can also define the map $\psi_s$ for $\tau_s$ and $\sigma_s$ 
in the same way the map $\psi$ was defined for $\sigma$ and $\tau$. 
However, notice that: 
\begin{lemma} The right hand side of Equation (\ref{Eqpsi}) only depends
on the strongly inner class $[\tau]$ of $\tau$.
\end{lemma}
\noindent{\bf Proof.} Since for any other real form $\tau_s$ strongly 
inner to $\tau$, the map $\sigma_s$ is given by $s.\sigma$, an orbit
$K.x$ will be mapped to $K.\sigma(x)$ by both $\sigma$ and $\sigma_s$.
\hfill {\bf Q.E.D.} \medskip

Let ${\mathcal I}_{\tau}$ be a finite set indexing the connected components 
of $Q$. It follows from our discussion that each element of this set defines a 
conjugacy class of involutions strongly inner to $\tau$. It is often
convenient to use elements of $K$ to represent elements of 
${\mathcal I}_{\tau}$. 

\begin{theorem} If $K$ acts freely on $\mu^{-1}(0)$, then $(X//K)^{\sigma}$ 
is diffeomorphic to the disjoint union of all $X^{\sigma_s}//K^{\tau_s}$ for
$s\in{\mathcal I}_{\tau}$ via the maps $\psi_s$. \label{th:main} 
\end{theorem}
\noindent{\bf Proof.} Consider a point $y\in (X//K)^{\sigma}$. First, 
we need to show that there exists a unique element $s\in{\mathcal I}_{\tau}$
such that $y$ is in the image of $\psi_s$. Let $y$ be represented by a
$\sigma$-stable $K$-orbit $K.x$ in $\mu^{-1}(0)$. If we let $\sigma(x)=s.x$, 
then the fact that $\sigma$ is involutive translates to $\tau(s)=s^{-1}$, 
so $s\in Q$. If we have a $\sigma$-fixed point $k.x$ in the orbit, then
we would have: 
$$k.x=\sigma(k.x)=\tau(k).\sigma(x)=\tau(k)s.x\ \ ,$$ which by the freedom 
of action would imply that $s=\tau(k^{-1})k$ and that $s$ is actually in $Q_0$.
If we do not have a $\sigma$-fixed point in the orbit, then $s$ is not 
in $Q_0$, but in a different connected component of $Q$. Then it is easy 
to check that the involutions $\tau_s$ and $\sigma_s$ defined by $s$ 
are such that the point $s.x$ is fixed by $\sigma_s$:
$$\sigma_s(s.x)=s.\sigma(s.x)=s\tau(s).\sigma(x)=s.x\ \ . $$
Therefore, $y$ is in the image of $\psi_s$. The statement about 
diffeomorphism follows now from Lemma \ref{leminj} together with a simple 
computation of the injectivity of $\psi_*$ on the level of tangent spaces. 
Finally, the fact that image of $\psi_s$ is both open and closed in 
$(X//K)^{\sigma}$ follows from \cite[Corollary 7.2]{OSS}.  \hfill {\bf Q.E.D.} 
\medskip

\begin{corollary} If $Q$ is connected, and $K$ acts freely on $\mu^{-1}(0)$, 
then the map $\psi$ is a diffeomorphism.  
\end{corollary}

\noindent{\bf Example.} Let us consider the case when $K=T$ is a torus. 
In this case, we can split $T=T_+\times T_-$ in such a way that 
$\tau$ fixed $T_+$ pointwise and acts as the inverse map on $T_-$. If the real 
dimension of $T_+$ is $n$, then $Q$ has exactly $2^n$ connected components.  
In this setting, a particular case of the above Corollary, when $T=T_-$ 
was proven by Goldin and Holm in \cite[Proposition 4.3]{GH}. 
\medskip 

\noindent{\bf Example.} Let $K={\rm SU}(n)$ and $\tau$ act as the complex 
conjugation. Since every unitary symmetric matrix $A$ can be represented as 
$\exp(\i B)$ for a real symmetric matrix $B$ \cite{Gant}, the set $Q$ 
is connected, and the map $\psi$ is bijective. 
\medskip

\noindent{\bf Example.} Let $K$ be simple of adjoint type (such as 
${\rm PU}(n)$). Then, as we remarked earlier, our notion of strong inner 
involution is the same as the standard notion of inner involution. The 
connected components of $Q$ correspond to equivalence classes
of symmetric spaces of $K$ defined by involutions inner to $\tau$.
The classification of symmetric spaces is well-known
\cite{Helga}, as well as the fact that inner classes of involutions are
in bijective correspondence with the order two automorphisms of the 
Dynkin diagram of $K$ \cite{AV}. For example, if $K={\rm PU}(n)$, and 
$\tau={\rm Id}$, there are $[n/2]+1$ connected components: a point and
$[n/2]$ projectivized grassmannians.    
\medskip

It is easy to verify that the connected components of $Q$ are parameterized
by the set $H^1(\Gamma, K)$, where $\Gamma\simeq{\mathbb Z}/2$ is generated by 
$\tau$. Let $Z$ be the center of $K$, so we have the exact sequence
$1\to Z\to K\to K/Z\to 1$, which induces a long exact sequence in cohomology
sets. Since, in principle, we know the answer for the simple adjoint case
and the torus case (see examples above), as well as we can easily compute 
$H^i(\Gamma, A)$, where $A$ is a compact abelian Lie group acted upon by 
$\Gamma$, we can use this sequence to compute $H^1(\Gamma, K)$.  

\medskip\noindent{\bf Example.} Let $\tau$ be inner to the identity automorphism. 
In this case we can assume that $\Gamma$ acts trivially on $K$ and therefore we 
can apply Theorem 6 from \cite[Chapter III]{Serre} to see that $H^1(\Gamma, K)$ 
can be identified with $T_2/W$, where $T_2$ is the set of elements of order $2$ 
in a maximal torus $T$ and $W$ is the Weyl group. For example, when 
$K={\rm SU}(n)$ and $\tau={\rm Id}$, there are exactly $[n/2]+1$ elements.  
\medskip

Let us now discuss the situations when the group $K$ does not act 
freely on $\mu^{-1}(0)$. In general, we cannot expect the
map $\psi$ to be injective. An easy example would be when $K={\rm SU}(2)$, 
$\tau$ the complex conjugation, $X=S^2\times S^2\times S^2$ with the 
diagonal action of $K$ and symplectic form (for example) 
$4\pi_1^*\omega+3\pi_2^*\omega+2\pi_3^*\omega$, where $\pi_i$ is the projection
onto the $i$-th factor and $\omega$ the standard invariant symplectic form on 
$S^2$. The anti-symplectic involution $\sigma$ on $X$ is simply the reflection
about the equatorial plane on each factor.  
The involutively reduced space in this case consists of two points and
the symplectically reduced space consists of a single point. The reason is that
each point on $X$ has a non-trivial stabilizer - the center of $K$. This 
example illustrates the fact that there are two equivalence classes of
triangles of side lengths $2,3,4$ under the action of the motion group
in ${\mathbb R}^2$ but only one such in ${\mathbb R}^3$. 
However, it is important to notice that if we pass to the group ${\rm PU}(2)$,
then the injectivity will hold, because ${\rm PU}(2)^\tau$ has now two
connected components. 

In general the following is true \cite[Proposition 2.3(iii)]{OSS}:
\begin{lemma} For each point $y\in X^{\sigma}//K^{\tau}$ the fiber 
$\psi^{-1}(\psi(y))$ is finite. 
\end{lemma}   

We see that the question of counting the number of points in the fiber
heavily depends on the stabilizer of the point. 
Let two points $x_1, x_2\in X^{\sigma}$ from 
different $K^{\tau}$ orbits be in the same $K$-orbit: $x_2=k.x_1$\ , 
$k\in K\setminus K^{\tau}$. This would imply that $k^{-1}\tau(k)\in K_1$ - 
the stabilizer of $x_1$ in $K$. Let us represent $k=k_0p$ with 
$p\notin K_1$, $p\in Q_0$, and $k_0\in K^{\tau}$. 
We immediately see that $p^2\in K_1$. So the problem of counting the number of
points in a given fiber amounts to counting the equivalence classes of elemens
from $Q_0$ satisfying the above conditions. 

Let us now turn to surjectivity questions. Let us assume that the momentum map 
$\mu$ is proper. According to \cite{SL}, the symplectic
quotient $X//K$ is stratified according to the conjugacy classes of stabilizers of
points in $\mu^{-1}(0)$. There is a unique open dense stratum called the 
{\it principal} stratum denoted by $(X//K)_{prin}$. The involution $\sigma$ on 
$X//K$ clearly maps a stratum to a stratum and thus 
$(X//K)_{prin}$ is $\sigma$-stable. Let us denote its fixed point set there by 
$(X//K)^{\sigma}_{prin}$. If we can show surjectivity for a (connected component) of 
$(X//K)^{\sigma}_{prin}$, then the surjectivity would also be valid for its closure 
in classical topology, because each connected component of 
$\mu^{-1}(0)^{\sigma}$ is compact. We shall assume that there is a point in the
principal stratum $\mu^{-1}(0)_{prin}$, which is $\sigma$-stable. This is a 
rather reasonable assumption, because otherwise there would be no points in 
$(X^{\sigma}//K^{\tau})_{prin}$. The stabilizer of $y$, which we call $H$, 
is therefore a $\tau$-stable subgroup of $K$. Let $N_H$ stand for the normalizer of 
$H$ in $K$ and let $L=N_H/H$ be the quotient group. It is easy to check that 
the $N_H$ is also $\tau$-stable and that the involution $\tau$ descends to the
quotient group $L$. 

Let us recall another result of Sjammar and Lerman \cite[Theorem 3.5]{SL}. 
They have established that the principal stratum $(X//K)_{prin}$ can be 
realized as the smooth Marsden-Weinstein reduced space of the set of 
points $X_H$ on $X$ with stabilizer $H$ by the action of the group $L$. 
Using this result, we can convert the questions of surjectivity for singular 
reduction to the smooth reduction which we dealt with in Theorem \ref{th:main}.  

\section{Example: flag manifolds.} 

Let $G$ be a connected complex semisimple Lie group of adjoint type 
with Lie algebra $\fg$. Let $\fh$ be a Cartan subalgebra of $\fg$ and 
let $\Delta$ be the corresponding root system.  
Let us fix a choice of positive roots $\Delta^+$ and let $\Sigma$ be the 
basis of simple roots. For any arbitrary subset $S\subset \Sigma$, we get 
a parabolic subgroup $P_S$ of $G$ in a standard way as follows. 
Every root $\beta \in \Delta$ has a unique decomposition 
$$ 
\beta=\sum_{\alpha\in\Sigma} n_\alpha(\beta)\alpha,
$$
where $n_{\alpha}(\beta)$ is a collection of either non-positive integers, 
in which case $\beta$ is a negative root, or non-negative integers, in 
which case $\beta\in\Delta^+$. Let $\Delta_S\subset\Delta$ stand for the 
set of roots which only involve simple roots from $S$ in the above 
decomposition. Let $\Delta_S^+$ be the subset of $\Delta_S$ consisting of 
positive roots and let $\overline{\Delta_S^+}$
be the complement of $\Delta_S^+$ in $\Delta^+$. Let us further define 
$$
\fp_S=\fh+\sum_{\alpha\in \Delta_S}\fg_\alpha + 
\sum_{\alpha\in \overline{\Delta_S^+}}\fg_\alpha,
$$
where $\fg_\alpha$ is the root space corresponding to $\alpha$. The first 
two summands in the above formula form the Levi factor of $\fp_S$ and the 
last one is the nilradical. Let $P_S$ stand for the parabolic subgroup of 
$G$ corresponding to $\fp_S$. If $P$ is any such subgroup then $X=G/P$ is a 
complex flag manifold. 

Let $B$ be the Killing form on $\fg$ and let $E_\alpha\in\fg_\alpha$ be
chosen such that 
$$
\begin{array}{cccc}
[E_\alpha, E_\beta] & = & m_{\alpha, \beta}E_{\alpha+\beta} & \mbox{ if }
\alpha+\beta\in \Delta \\
{} & = & H_\alpha & \mbox{ if } \alpha=-\beta \\
{} & = & 0 & \mbox{otherwise}    
\end{array},
$$
where $H_\alpha$ is the unique element of $\fh$ defined by $B(H,
H_\alpha)=\alpha(H)$ for all $H\in\fh$. In addition we require that the
constants $m_{\alpha, \beta}$ be real and satisfy $m_{-\alpha,
-\beta}=-m_{\alpha, \beta}$. We take the compact real form 
$\fk$ of $\fg$ as the span of $\i H_\alpha$, $E_\alpha-E_{-\alpha}$, and
$\i(E_{\alpha}+E_{-\alpha})$. Let $\theta$ stand for the corresponding Cartan
involution. Let $\fg=\fk\oplus\fq$ be the Cartan decomposition 
of $\fg$ into the $\pm 1$ eigenspaces of the involution $\theta$. Let also $\fg=
\fk\oplus\fa\oplus\fn$ be the Iwasawa decomposition corresponding to our choices
of $\Delta^+$ and $\fk$. Let $G=KAN$ be the corresponding decomposition on the 
group level. 

Each isomoprhism class of real forms of $G$ is represented by a Satake diagram 
$D$, as explained e.g. in \cite{Helga}. A Satake diagram is the Dynkin diagram
for $\fg$ such that some vertices are painted black and certain white
vertices are paired by arrows. As explained in \cite{FothLu}, we can 
construct an involution $\tau$ on $G$ such that $G^\tau$ is a specific 
representative of a real form corresponding to $D$, called an Iwasawa real form, 
which has the following properties: 

\noindent(1) $\tau$ commutes with $\theta$: $\tau\circ\theta=\theta\circ\tau$.
Then $\fg^\tau=\fk^\tau + \fq^\tau$. 

\noindent(2) $\fh^{\tau}$ is a maximally non-compact Cartan sub-algebra in $\fg^\tau$.

\noindent(3) If we denote by $N^{\tau}$ the subgroup of $N$ consisting of
elements fixed by $\tau$, then $$G^\tau=K^{\tau}A^{\tau}N^{\tau}$$ is an Iwasawa
decomposition of $G^{\tau}$. 

Note that the group $N$ is not stabilized by $\tau$ in general, but only when 
$D$ corresponds to the so-called {\it quasi-split} real form \cite{AV}. 

Since we are working with an adjoint group, the notion of strongly innner 
involution is the same as inner, and we will use another result of \cite{AV} 
which asserts that there is a unique, up to conjugacy, quasi-split real form  
represented by a Satake diagram $D_{qs}$ with no black vertices for each inner 
class of real forms. Let us denote the corresponding Iwasawa involution 
$\tau_{qs}$. We may assume that $\tau_{qs}$ commutes with $\tau$.  
     
\medskip\noindent{\bf Example.} Let us give examples of Iwasawa real forms
for the $A_{n-1}$ case, when $G={\rm PGL}(n, {\mathbb C})$. There are two inner 
classes of real forms, I and II. The class I contains real forms isomorphic
to ${\rm PU}(p,q)$'s for $p+q=n$, $p\le q$, and the class II contains the 
split real form ${\rm PGL}(n, {\mathbb R})$ and, when $n=2m$ is even, 
a real form isomorphic to ${\rm PGL}(m, {\mathbb H})$. Let us write down 
specific Iwasawa involutions. Let 
$$
Q_p=\left( \begin{array}{ccc} 
0_{p\times p} & 0_{p\times(q-p)} & Y_{p} \\
0_{(q-p)\times p} & 1_{q-p} & 0_{(q-p)\times p} \\
Y_{p} & 0_{p\times (q-p)} & 0_{p\times p}
\end{array}\right)\ \ ,
$$
where $Y_p$ is the $p\times p$ matrix with ones on the anti-diagonal.   
Then the involution $$ \tau_p(A)= Q_p(\overline{{}^tA^{-1}})Q_p $$
determines a real form isomorphic to ${\rm PU}(p,q)$. In this inner class the 
split real form is given by $\displaystyle{{\rm PU}\left( \left[ \frac{n}{2}\right], 
\left[ \frac{n+1}{2}\right]\right)}$.   

The class II contains the split real form ${\rm PGL}(n, {\mathbb R})$ 
with the corresponding Iwasawa involution being the complex conjugation and, 
when $n=2m$ is even, we have the real form isomorphic to ${\rm PGL}(m, {\mathbb H})$ 
with the corresponding Iwasawa involution given by 
$$ \tau_{\mathbb H}(A) = Q_{\mathbb H} {\bar A} Q_{\mathbb H}^{-1}\ \ , $$
where $Q_{\mathbb H}$ is $2m\times 2m$ block-diagonal matrix: 
$$Q_{\mathbb H}={\rm diag}(J, -J, J, ..., (-1)^m J), \ \ \ 
J=\left( \begin{array}{cc} 0 & 1 \\ -1 & 0 \end{array} \right)\ \ .$$
\medskip

Now let us recall some classical results about the action of the real
form $G^\tau$ on the complex flag manifold $X=G/P$, which are due to Joseph
Wolf \cite{wolf1}. As it is standard, we view points in $X$ as a $G$-conjugates
of $\fp$ via the correspondence $gP\leftrightarrow \text{Ad}(g)\fp$.
We denote by $P_x$ (resp. by $\fp_x$) the corresponding parabolic subgroup of $G$ 
(resp. the parabolic subalgebra of $\fg$). Wolf showed that
there exists a $\tau$-stable Cartan subalgebra $\fh_x\subset\fp_x$ of
$\fg$, a positive root system $\Delta^+_x$ compatible with $\fh_x$, and a set
$S_x$ of simple roots, such that $\fp_x=\fp_{S_x}$ and $P_x=P_{S_x}$. 
The real co-dimension of the $G^{\tau}$-orbit $G^{\tau}(x)$ through $x$ is 
equal to the cardinality of the intersection of $\overline{\Delta^+_{S_x}}$ with 
$\tau\overline{\Delta^+_{S_x}}$. The number of $G^{\tau}$-orbits on $X$ is finite, 
and there is a unique closed orbit $X_0$, which is contained in the closure of 
every $G^{\tau}$-orbit. If $G^{\tau}(x)$ is the closed orbit, then there is an 
Iwasawa decomposition $G^{\tau}=K^{\tau}A^{\tau}N^{\tau}$ such that $G^{\tau}\cap P_x$ 
contains $H^{\tau}N^\tau$, whenever $H^{\tau}$ is a Cartan subgroup of $G^{\tau}$ 
containing $A^{\tau}$. Moreover, if $K_0$ is any maximal compact subgroup of $G^{\tau}$, 
then $K_0$ is transitive on $X_0$. 

It will be of a particular interest to us to investigate those cases when the 
real dimension of $X_0$ is half the real dimension of $G/P$. It was proved in 
\cite{wolf1} that the closed $G^{\tau}$-orbit $X_0=G^{\tau}(x)$ satisfies
${\rm dim}_{\mathbb R}(X_0)={\rm dim}_{\mathbb C}X$ if and only if the following
equivalent conditions hold: 

\medskip\noindent(1) $\tau\overline{\Delta^+_{S_x}}=\overline{\Delta^+_{S_x}}$.

\noindent(2) $\tau(\fp_x)=\fp_x$, $\tau(P_x)=P_x$. 

\noindent(3) $X$ a projective variety defined over ${\mathbb R}$ and 
$X_0$ is the set of real points.
\medskip

It is therefore appropriate in such situations to refer to $X_0$
as a {\it real flag manifold}. In the case when $P=B$, the Borel subgroup, 
the condition ${\rm dim}_{\mathbb R}(X_0)={\rm dim}_{\mathbb C}X$
holds if and only if the Satake diagram of the symmetric space $G^{\tau}/K^{\tau}$
contains no painted vertices and $\fk^{\tau}$ does not contain a simple ideal of
$\fg^{\tau}$. 

Let $X_0$ be the closed $G^{\tau}$-orbit on $G/P$. We can assume that $P$ is chosen 
in such a way that $X_0$ is the orbit through the base point of $G/P$. If 
${\rm dim}_{\mathbb R}X_0={\rm dim}_{\mathbb C}X$, then the previous discussion implies 
that $\tau(\fp)=\fp$ and that $\tau\overline{\Delta^+_S}=\overline{\Delta^+_S}$. 
According to \cite{MacDo}, the set of simple roots $\Sigma$ decomposes into the
disjoint union of two subsets $\Sigma_0$ and $\Sigma_1$ such that for any 
$\alpha\in\Sigma_0$ we have $\tau(\alpha)=-\alpha$ and for any
$\alpha\in\Sigma_1$ we have 
$$\tau(\alpha)= \mu(\alpha)+\sum_{\beta\in\Sigma_0}c_{\alpha, \beta}\beta,$$
where $\mu$ is an involution on the set $\Sigma_1$ and $c_{\alpha, \beta}$ are 
non-negative integers. Let $S$, as before, be the subset of simple roots defining
the parabolic subgroup $P$. First of all, the condition 
$\tau\overline{\Delta^+_S}=\overline{\Delta^+_S}$ implies that 
$\Sigma_0\subset S$. Let us denote $S_1:=S\cap \Sigma_1$. For any $\alpha\in S$ 
let us define by $h_\alpha\in\fa$ the unique element with the property that 
for any $\beta\in S$ we have $\beta(h_\alpha)=\delta_{\alpha, \beta}$. 
We define $$\tilde{\lambda} = \sum_{\alpha\in \Sigma_1\setminus S_1}h_{\alpha} \ \ 
\in\fa.$$ One checks straightforwadly that:
\begin{lemma} $$\tau(\tilde{\lambda})=\tilde{\lambda}$$
\end{lemma}

Let us now identify $\i\fk$ with $\fk^*$ using the imaginary part of the
Killing form and let us denote by $\lambda$ the image of $\tilde{\lambda}$
under this identification. Due to the fact that $\tau$ is complex anti-linear, 
and this identification uses $\i$, we have $\tau(\lambda)=-\lambda$. 

Let us denote $K_P=K\cap P$, then we have an identification $G/P\simeq K/K_P$. 
The element $\lambda\in\fk^*$ defines a left-invariant one-form on $K$, which 
we also call $\lambda$. Let $\omega$ be the unique two-form on $X=K/K_P$ such
that $p^*\omega=-d\lambda$, where $p$ is the projection $K\to X$. It is a
standard fact that $\omega$ defines a symplectic structure on $X$. Since the
involution $\tau$ preserves the parabolic subgroup $P$, it descends to
an involution denoted by $\sigma$ on $X=G/P$ and thus we have: 

\begin{theorem} When $\emph{dim}_{\mathbb R}(X_0)=\emph{dim}_{\mathbb C}X$, there exists a
symplectic structure on $X$ and an anti-symplectic involution $\sigma$ of $X$ such
that its fixed point set $X^{\sigma}$ is the closed $G^{\tau}$-orbit $X_0$.  
\end{theorem} 
\noindent{\bf Proof.} 
The fact that $\tau(\lambda)=-\lambda$ immediately implies that the
involution $\sigma$ on $X$ is anti-symplectic and satisfies 
$\sigma(k\cdot x)= \tau(k).\sigma(x)$, where $k\in K$. The fixed point set of 
$\sigma$ is non-empty, because, for example, it contains the closed
$K^{\tau}$-orbit, which is $X_0$ by \cite{wolf1}. 
The rest follows from \cite[Example 2.9]{OSS}. 
\hfill {\bf Q.E.D.}\medskip

\noindent{\bf Remark.} O'Shea and Sjamaar studied real flag manifolds
as examples of fixed loci of anti-symplectic involutions on what they called 
{\it symmetric} co-adjoint orbits. We refer to \cite{OSS} for more detail. \medskip

Therefore a large class of examples for which we can apply our previous
results is given by a product of (partial) flag manifolds for the group $G$ such 
that each term in the product is equipped with an anti-symplectic involution
compatible with the real form $G^\tau$. 

In fact, for each isomorphism class of real forms represented by a Satake diagram
$D$ and an Iwasawa involution $\tau$, we can classify all standard parabolic 
subgroups, which are stabilized by $\tau$. The condition is that the corresponding
subset $S$ of real roots has to contain all black vertices from $D$ and, 
in addition, if one white vertex from a pair connected by an arrow is in $S$, 
then the other one should be in $S$ as well. So the classification boils down to
a simple combinatorics involving the Satake diagram. Moreover, if we take the 
quasi-split $D_{qs}$ in the same inner class together with a commuting Iwasawa
involution $\tau_{qs}$, then one can easily see that if a standard parabolic 
subgroup $P_S$ is stabilized by $\tau$, then it will also be stabilized by 
$\tau_{qs}$. 
 
For each such standard parabolic subroup $P_S$ stabilized by $\tau$, the flag
manifold $G/P_S$, according to our previous discussion, will have an anti-symplectic
involution $\sigma$ compatible with $\tau$ such that $(G/P_S)^{\sigma}$ is 
a non-empty lagrangian submanifold of $G/P_S$. It follows that the involution
$\sigma_{qs}=s\sigma$, determined by $\tau_{qs}={\rm Ad}_s\circ\tau$ will also have a 
non-empty lagrangian fixed point set on $G/P_S$. Therefore, after identifying
the flag manifolds of the form $G/P_S$ with the co-adjoint orbit carrying the same
symplectic form $\omega$, we arrive at the following result: 

\begin{proposition} Let ${\mathcal O_i}\subset\fk^*$, $1\le i\le k$ be a collection
of co-adjoint orbits such that ${\mathcal O}_i\simeq G/P_i$ and for each
$i$, the standard parabolic subgroup $P_i$ is $\tau$-stable. Then the
equation $A_1+\cdots + A_k=0$, $A_i\in {\mathcal O}_i^{\sigma_{qs}}$
has a solution if and only if the same equation has a solution with   
$A_i\in {\mathcal O}_i^{\sigma}$. 
\end{proposition}

\noindent{\bf Proof.} 
If $A_1+\cdots +A_k=0$ has a solution with 
$A_i\in {\mathcal O}_i^{\sigma}$, then our preceeding discussion implies
that ${\mathcal O}_i^{\sigma_{qs}}$ is non-empty. The rest follows from
\cite[Theorem 3.1(i)]{OSS}. \hfill {\bf Q.E.D.} \medskip

Thus we can always push the problem of solving the problem about the sum of 
elements with prescribed spectra for a given real form of the group to the 
same problem for the quasi-split real form in that inner class. 

\medskip\noindent{\bf Remark.} The inequalities for the additive problem 
are given in \cite{OSS}. The equivalence to the
multiplicative problem for the quasi-split real form
with prescribed singular values follows from \cite{AMW}. 
\medskip

\noindent{\bf Example.} Let us consider a complex flag manifold
$X={\rm Fl}_{\mathbb C}(2m_1, 2m_2, ..., 2m_k)$ which parameterizes complex flags 
$$V_{2m_1}\subset V_{2m_2}\subset\cdots\subset V_{2m_k}\simeq{\mathbb C}^{2m}\ \ ,$$
${\rm dim}_{\mathbb C} V_{2i} = 2i$ and $0< m_1 <\cdots < m_k =m$.   
We will use the natural embedding
$$\iota:\ X\hookrightarrow \prod_{i=1}^k {\rm Gr}_{\mathbb C}(2m_i, 2m)
$$ to get a symplectic structure on $X$, and the symplectic structure on
each of the grassmannians is a positive multiple of the one that
comes from the standard Pl\"ucker embedding
and the Fubini-Study form on the projective space of dimension $2m$ 
choose $2m_i$. 

We shall identify ${\mathbb C}^{2m}$ with the right quaternionic space 
${\mathbb H}^m$ as follows. 
The point $(z_1, ..., z_{2m})\in {\mathbb C}^{2m}$ corresponds to the point 
$(q_1, ..., q_m)\in{\mathbb H}^m$ 
if $q_i=z_{2i-1}+{\bf j}z_{2i}$ for $1\le i\le m$. Using this identification, let $J$ be 
the real operator on ${\mathbb C}^{2m}$ which comes from the right
multiplication by ${\bf j}$ on the space ${\mathbb H}^m$. Since $J^2=-\mbox{Id}$, the action of $J$ 
on ${\mathbb C}^{2m}$ extends to an involution $\sigma$ on all complex partial flag manifolds. 
However, this involution is only real and not a complex diffeomorphism. If all the 
weights of the subspaces are even, then the fixed point set of $\sigma$ is clearly 
the quaternionic (partial) flag manifold: 
$$
({\rm Fl}^{{\mathbb C}}(2m_1, ..., 2m_k))^{\sigma}={\rm Fl}^{{\mathbb H}}(m_1, ..., m_k).
$$
	Moreover, if $\omega$ is an invariant K\"ahler form on the complex flag 
manifold, then the fixed point set of $\sigma$ is a Lagrangian submanifold with 
respect to $\omega$. For example, when $X=Gr_{{\mathbb C}}(2,4)$, we have 
$X^\sigma={\mathbb H{\mathbb P}}^1\simeq {\rm S}^4$.  As usual, we view the group 
${\rm PGL}(m, {\mathbb H})$ as a real form of the complex semi-simple group 
${\rm PGL}(2m, {\mathbb C})$, with the corresponding Satake 
diagram \cite{Helga} having odd numbered vertices painted black and no arrows.  
The multiplicity of each restricted root is $4$. In terms of matrices, we have the embedding
\begin{equation}
\nu: {\rm PGL}(m, {\mathbb H})\hookrightarrow {\rm PGL}(2m, {\mathbb C}), \ 
A+B{\bf j}\ \mapsto \left(
\begin{array}{cc}
A & B \\ -{\bar B} & {\bar A}
\end{array}
\right).
\label{e:nu}
\end{equation}
If we let $J$ be the $2m\times 2m$ matrix 
$$
\left(
\begin{array}{cc}
0 & 1 \\ -1 & 0
\end{array}
\right) \ \  , $$
then one can define the involution $\tau$ on ${\rm PGL}(2m, {\mathbb C})$ by 
$\tau(C)=-J{\bar C}J$, which defines the real form ${\rm PGL}(m, {\mathbb H})$. 
(To simplify the discussion, we use a different $\tau$ than earlier.)
One can see that $\sigma$ and $\tau$ are compatible since both have the same origin 
and the action of ${\rm PGL}(m, {\mathbb H})$ on the fixed 
point set of $\sigma$ on flag manifolds comes from the action of ${\rm PGL}(2m, {\mathbb C})$.  

Now we will conjugate $\tau$ to the Iwasawa involution $\tau_{\mathbb H}$ defined
previously, and consider together with the involution $\tau_{qs}$ which is, simply, 
the complex conjugation. Any standard parabolic subgroup defined by any set $S$ 
consisting of all odd-numbered and an arbitrary subset of even-numbered vertices in the 
Satake diagram will be simultaneously $\tau_{\mathbb H}$- and $\tau_{qs}$-stable. 
Then, as usual, we can take $X$ to be a product of a number of such $G/P_S$ 
and arrive to a series of examples when $(X//K)^\sigma$ will have separate connected 
components corresponding to the involutive quotients of the weighted $k$-fold product
of quaternionic flag manifolds by the diagonal action of groups ${\rm PSp}(m)$ and 
weighted $k$-fold product of real flag manifolds by the diagonal action of ${\rm PO}(2m)$. 

\medskip \noindent {\bf Example.} Let us take $m=2$ and $m_1=1$, and $X$ the 
$m$-fold product of complex grassmannians of two-planes in ${\mathbb C}^4$:
$$ X= {\rm Gr}_{\mathbb C}(2,4)\times \cdots \times {\rm Gr}_{\mathbb C}(2,4). $$
Each factor ${\rm Gr}_{\mathbb C}(2,4)$ carries its own invariant symplectic form, 
a positive multiple of a standard one. The fixed point set of the involution $\sigma$ 
in this case on each factor is the quaternionic projective line 
${\mathbb H}{\mathbb P}^1\simeq {\rm S}^4$ and the quotient 
$$X^{\sigma}//K^{\tau}=\left( {\rm S}^4\times\cdots \times {\rm S}^4 \right) 
// {\rm PSp}(2)\ \ ,$$
considered in detail in \cite{FL}, is identified with the moduli space of polygons 
in ${\mathbb R}^5$. Indeed, a point in ${\rm S}^4$ gives a direction in ${\mathbb R}^5$, 
the weight on each factor is the corresponding side length, and 
${\rm PSp}(4)={\rm SO}(5)$. On the other hand, the real locus of 
$$\left( {\rm Gr}_{\mathbb C}(2,4)\times \cdots \times {\rm Gr}_{\mathbb C}(2,4)
\right) //{\rm SU}(4)$$ defined by $\sigma$ will have another connected component
corresponding to the involutive quotient 
$$\left( {\rm Gr}_{\mathbb R}(2,4)\times \cdots \times {\rm Gr}_{\mathbb R}(2,4)
\right) //{\rm PO}(4).$$

\section{Invariant theory quotients and involutions.}

Let $K$ be a compact connected Lie group and let $G$ be its complexification. 
Then $G$ is connected linearly reductive complex Lie group.
Let $\tau$ be an antiholomorphic involution on $G$ defining a real form 
$G^{\tau}$.  Then we have $K^{\tau}=K\cap G^{\tau}$.
Let $\rho: G\to {\rm GL}(N,{\mathbb C})$ be a rational 
representation and let $\theta$ be a Cartan involution of 
${\rm GL}(N, {\mathbb C})$ commuting with $\tau$. This condition is not 
very restrictive, since any anti-holomorphic involution on $G$ is conjugate 
to a one commuting with $\theta$. Let $V$ be a complex vector space
of dimension $N$ on which ${\rm GL}(N, {\mathbb C})$ acts in the usual way,
equipped with an anti-holomorphic involution $\sigma$ compatible with $\tau$.
We can always assume that the image of $K$ under $\rho$ is contained in 
${\rm U}(N)={\rm GL}(N, {\mathbb C})^\theta$ and thus $K$ acts by self-adjoint
operators.   

Let us recall some of the results from Richardson-Slodowy \cite{RS}. 
Denote by ${\mathcal M}$ the set of minimal vectors in $V^{\sigma}$
for the $G^{\tau}$-action
and  consider the map $\pi: V^{\sigma}\to V^{\sigma}//G^{\tau}$ which is the 
Luna quotient of $X^{\sigma}$ by $G^{\tau}$ whose points are the closed orbits 
equipped with the quotient topology. In {\it loc.cit.} it was shown that 
the inclusion ${\mathcal M}\hookrightarrow V^{\sigma}$ induces a homeomorphism 
of the orbit space ${\mathcal M}/K^{\tau}$ with the quotient $V^{\sigma}//G^{\tau}$. 
The set ${\mathcal M}$ coincides with the $\sigma$-invariant minimal
vectors for the $G$-action on $V$. For any closed $G^{\tau}$-stable real-algebraic 
subset $Z\subset V^{\sigma}$ the induced map $Z/K^{\tau}\to Z//G^{\tau}$ is a homotopy 
equivalence. Let $\beta: V^{\sigma}//G^{\tau}\to V//G$ be the natural map. For each 
$\xi\in V^{\sigma}//G^{\tau}$, the fiber $\beta^{-1}(\beta(\xi))$ is finite.  
 
Let now $\omega$ stand for the Fubini-Study form on ${\mathbb C{\mathbb P}}^{N-1}$
for which $\sigma$ on ${\mathbb C{\mathbb P}}^{N-1}$ satisfies 
$\sigma^*{\omega}=-{\omega}$. Let $X\hookrightarrow {\mathbb C{\mathbb P}}^{N-1}$
be a smooth $\sigma$-stable projective variety equipped with a linear 
action of $G$ via $\rho$. 

By standard results in the invariant theory \cite{GIT}, there is a
momentum map $\mu: X\to \fk^*$ such that the inclusion of $\mu^{-1}(0)$ 
into the subset $X^{ss}$ of semi-stable points induces a 
homeomorphism:
\begin{equation} \phi:\ \mu^{-1}(0)/K\to X//G\ \ \ , \label{eqgithomeo} 
\end{equation} 
where the right hand side is the categorical quotient of $X$ by the 
action of $G$. According to our previous discussion, the involution
$\sigma$ leaves $\mu^{-1}(0)$ invariant and descends to the left-hand
side of (\ref{eqgithomeo}). It is also easy to check that the 
involution $\sigma$ descends to the categorical quotient $X//G$, 
and the above homeomorphism respects these involutions.  

According to \cite{RS}, we have a homeomorphism 
$$ \gamma:  
\left( X^{\sigma}\cap \mu^{-1}(0)\right) /K^{\tau}\to X^{\sigma}//G^{\tau}\ \ .$$
Thus we have the following commutative diagram of continuous maps:
$$\begin{array}{ccc}
X^{\sigma}//K^{\tau} & \stackrel{\psi}{\longrightarrow} & (X//K)^{\sigma} \\
\gamma \downarrow & & \phi \downarrow \\
X^{\sigma}//G^{\tau} & \stackrel{\eta}{\longrightarrow} & (X//G)^{\sigma}
\end{array}\ \ ,$$
where $\gamma$ and $\phi$ are homeomorphisms. In particular, each of the four
spaces in the diagram is homeomorphic to a closed semi-algebraic set. 

\medskip\noindent{\bf Example.} Let us return to the example of
quaternionic flag manifolds from Section 3. 
Let $m_i<m$ be two positive integers and let 
$\displaystyle{k_i=\left( {\begin{array}{c} 2m \\ 2m_i \end{array}}\right)}$.
Consider the Pl\"ucker embedding 
${\rm Gr}_{\mathbb C}(2m_i, 2m)\hookrightarrow {\mathbb C{\mathbb P}}^{k_i-1}$. 
It is easy to see that the action defined by $J$ on ${\mathbb C}^{2m}$ 
lifts to a complex-conjugate involution denoted by $\sigma$ on 
${\mathbb C}^{k_i}$ for the reason that $2m_i$ is an even number. 
The compatibility of actions is straightforward. In particular, 
we can think of the moduli space of polygons in ${\mathbb R}^5$ as
the real invariant quotient of $({\mathbb H{\mathbb P}}^1)^n$ by the diagonal 
action of the group ${\rm PGL}(2, {\mathbb H})$. \medskip
  
Returning to the general case, we can translate the results of our 
previous discussion to the map $\eta$ and for each involution 
$\tau_s={\rm Ad}_s\circ\tau$ strongly inner to $\tau$ construct the 
corresponding involution $\sigma_s$ and map $\psi_s$ as in Section 2. 
Note that the element $s$ can be chosen from $K$. 

\begin{theorem} If $X//G$ is a smooth variety, then 
the manifold $(X//G)^{\sigma}$ is homeomorphic to the disjoint union of 
$\psi_s(X^{\sigma_s}//G^{\tau_s})$, where $s$ runs through the 
set ${\mathcal I}_{\tau}$ as in Section 2. If, moreover, $G$ acts freely on 
stable points, then  each $\psi_s$ is injective.  
\end{theorem}

When $X//G$ has quotient singularities, then according to our results in 
Section 2, the maps $\psi_s$ still cover $(X//G)^{\sigma}$, although the
(finite) fibers over the singular points may differ from fibers
over the smooth points.  

\section*{Acknowledgements}

I am grateful to Jiang-Hua Lu for helpful comments and   
suggestions. I would also like to thank Paul Bressler, Sam Evens, 
Dipendra Prasad, and Reyer Sjamaar for useful discussions and correspondence. 
I acknowledge the NSF grant DMS-0072520.


\end{document}